\documentclass[journal]{IEEEtran}
\ifCLASSINFOpdf
\else
\fi
\usepackage[comma,numbers,square,sort&compress]{natbib}
\usepackage{algorithm} 
\usepackage{algorithmic} 
\usepackage{multirow} 
\usepackage{xcolor}
\usepackage{amsmath}
\newtheorem{theorem}{\textbf{Theorem}}
\newtheorem{lemma}{\textbf{Lemma}}
\newtheorem{corollary}{\textbf{Corollary}}
\newtheorem{remark}{\textbf{Remark}}
\newtheorem{definition}{\textbf{Definition}}

\usepackage{multirow}
\usepackage{url}
\usepackage{subfigure}
\usepackage{fancyhdr}
\usepackage{amsmath}
\usepackage{multirow}
\usepackage{amssymb }
\usepackage{color}
\usepackage{graphics} 
\usepackage{graphicx}

\begin{document}
%
\title{\LARGE \bf
Minimal Structural Perturbations for Network Controllability: Complexity Analysis
}
%
%
%

\author{Yuan Zhang and Tong Zhou$^{\dag}$
\thanks{*This work was supported in part by the NNSFC under Grant 61573209 and 61733008. {\bf This work is just an extension of the conference paper \cite{Y_Zhang_2017}}.}
\thanks{$^{\dag}$Yuan Zhang and Tong Zhou are with the Department of Automation and TNList, Tsinghua University, Beijing, 100084, P.~R.~China
        {(email: {\tt\small zhangyuan14@mails.tsinghua.edu.cn, tzhou@mail.tsinghua.edu.cn}).}}%
}
\maketitle

\begin{abstract} Link (edge) addition/deletion or sensor/actuator failures are common structural perturbations for real network systems. This paper is related to the computation complexity of minimal (cost) link  insertion, deletion and vertex deletion with respect to structural controllability of networks. Formally,  given a structured system, we prove that: i) it is NP-hard to add the minimal cost of links (including links between state variables and from inputs to state variables) from a given set of links to make the system structurally controllable, even with identical link costs or a prescribed input topology; ii) it is NP-hard to determine the minimal cost of links whose deletion deteriorates structural
controllability of the system, even with identical link costs or when the removable links are restricted in input links. It is also proven that determining the minimal cost of inputs whose deletion causes structural uncontrollability is NP-hard in the strong sense. The reductions in their proofs are technically independent.  These results may serve an answer
to the general hardness of optimally designing (modifying) a structurally controllable
network topology and of measuring controllability robustness against link/actuator
failures. Some fundamental approximation results for these related problems are also provided.  
\end{abstract}

\section{Introduction}
Recently, the design of large scale systems has attracted  much interest with the emergence of complex networks, such as power networks, biological transduction networks \cite{nature 2011}, gene regulation networks \cite{Xiongjie_2014}, etc. One fundamental objective is to design a network that ensures controllability and observability \cite{A_Olsehvsky_2014}, \cite{T_H_Summers_2016},  \cite{zhou_2015}, \cite{output-degree}, \cite{Y_Zhang_2016}.
  Among the related problems, the input selection problem has received much attention in \cite{A_Olsehvsky_2014}, \cite{T_H_Summers_2016}, \cite{Zhou_minimal_control_2016}, \cite{C_Commault_2013}, \cite{nature 2011}, \cite{Sergio_Pequito_2016}, \cite{Sergio_Pequito_2016_cost}. Specially, it is known that determining the minimal actuated states to ensure controllability for a numerical system is NP-hard \cite{A_Olsehvsky_2014}. However, if we ignore the exact parameters of the system matrices and only focus on their zero-nonzero sparsity patterns, then the same problem of determining the minimum number of actuated states to ensure structural controllability can be done in polynomial time \cite{Sergio_Pequito_2016} using some graph theoretical operations.  Apart from the binary concept of controllability, researchers also develop some heuristic methods to select inputs to optimize certain control energy related metrics \cite{T_H_Summers_2015}, \cite{T_H_Summers_2016}.
Compared to the abundant research on input selection problems, less attention has been paid to the design of the (autonomous) network topology. In this paper, we are interested in the following
questions which often emerge in designing the topologies of
networks:
given a system, 1) if the system is uncontrollable, how to adjust links between state variables, or from the existing inputs to state variables rather than adding extra inputs, to make the system controllable? 2) inverse to 1), if the system is controllable, how to identity the subsets of links/actuators whose removal would destroy the system controllability?
Here links could correspond to interacting connections, communication channels, connectivity paths etc. in practical networks, such as multi-agent systems, complex communications networks, transportation systems \cite{multivariable_graph_1988}. Some preliminary work concerning Problem 1) can be found in \cite{zhou_2015}, where it illustrates how to transform a specific uncontrollable networked system to be a controllable one (in numerical sense) by adjusting subsystem connections, yet systematic methods to do this transformation need further study.
As an inverse problem of Problem 1), Problem 2) can provide information concerning the robustness of system topologies, or the 'Achilles heel' link/actuator sets, i.e., elements whose absence will make the system uncontrollable.  A simple classification for network links can be found in \cite{nature 2011} according to the effects of their absence on the number of driver nodes needed to ensure controllability.  Since controllability and observability are closely related security of cyber-physical systems, Problem 2) is significant to determine whether a system is resilient under malicious link/actuator attacks with bounded cardinality.


Structural controllability is only related to the zero-nonzero patterns of the associated system matrices \cite{Lin_CT_1974}, which serves as an alternative notion for controllability if we have no access to the exact value of the link weights of the networks.  
The problem of modifying a network by adding links between state vertices to make the network controllable by one single input has been considered in \cite{wangwenxu}. The problem of building a structurally observable system with minimum link cost and robustness consideration has been studied in \cite{S_Pequito_2017} under the assumption that all state variables have zero-cost self-loops. 
 Robustness of controllability and observability under structural disturbances have been discussed in \cite{Robustness_IJC}. \cite{C_Comault_Observability} considers observability preservation under sensor failure; later \cite{M_A_Rahimian_failures_2013} studies controllability preservation under simultaneous failures in both the communication links and the agents. These works mainly focus on classification of links and agents according to the influence of their failures on observability or controllability. However, computation complexity concerning on the associated optimization problems, to the best of our knowledge, has not been formally established in literature.

In this paper, we study the computation complexity of the optimization versions of the link (edge) insertion/deletion and
actuator deletion subject to structural controllability. These problems are significant to understanding the `distance' between structural controllability and structural uncontrollability \cite{R_Eising_1984}, \cite{Y_Zhang_2017}.  Since these problems are combinatorial problems at first look, understanding their computation complexity is important.   Our main contributions are three complexity results concerning the minimal (cost) structural perturbations for network controllability. To be specific, given a structured system, we prove that: i) it is NP-hard to add the
minimal cost of links (including links between state variables
and from inputs to state variables) from a given set of links to
make the system structurally controllable, even with identical
link costs or a prescribed input topology; ii)  it is NP-hard to
determine the minimal cost of links whose deletion deteriorates
structural controllability of the system, even with identical link
costs or when the removable links are restricted in input links; and iii) it is NP-hard in the strong sense to determine the minimal cost of actuators whose deletion causes uncontrollability.  While these three problems are conceptually related, the proofs in their reductions are technically independent. The first result is in sharp contrast to the recently known fact that selecting the minimal number (cost) of states to be actuated to ensure structural controllability can be solved in polynomial time \cite{Sergio_Pequito_2016}. The second result means that, it is impossible to determine the controllability `robustness' against link failures in polynomial time under the common conjecture $P\neq NP$. Strong NP-hardness means that, there is no quasi-polynomial time algorithms for the third problem unless $P=NP$.  Some fundamental approximation results for these related problems are also provided. For example, we show that a $2$-approximation polynomial time algorithm exists for the first problem, and the second problem has the same multiplicative approximation factor as that of the minimal cost $1$-blocker problem.  These results may serve an answer to the general hardness of optimally designing (modifying) a structurally
controllable network topology and of measuring controllability robustness against link/actuator failures.

The rest of this paper is organized as follows. Section II provides some preliminaries and introduces the problems studied in this paper. Sections III, IV and V respectively give the intractability and approximation results for the associated link insertion, link deletion and actuator deletion problems respectively.  The concluding remarks are included in Section VI.

\section{Preliminaries and Problem Formulation}

\subsection{Concepts in Graph Theory}
Given a digraph $\mathcal{G}=(\mathcal{V},\mathcal{E})$, a path from $v_i$ to $v_j$ is a sequence of edges $\{(v_i,v_{i+1}),(v_{i+1},v_{i+2}),...,(v_{j-1},v_j)\}$ without repeated vertices. A digraph is said to be strongly connected, if for any two vertices $v$ and $w$ of this digraph, there is a path from $v$ to $w$ and from $w$ to $v$, i.e., $v$ and $w$ can be reachable from each other. A strongly connected component (SCC) of $\mathcal{G}$ is a subgraph of $\mathcal{G}$ that is strongly connected and is maximal in the sense that no additional edges or vertices from $\mathcal{G}$ can be included in the subgraph
without breaking its property of being strongly connected. 

For a graph $\mathcal{G}$, $V(\mathcal{G})$ denotes the vertex set of graph $\mathcal{G}$, $E(\mathcal{G})$ the edge set. Given a digraph $\mathcal{G}=(\mathcal{V},\mathcal{E})$ with edge costs (weights) $c: \mathcal{E}\rightarrow \{0\}\cup \mathbb{R}^+$, denote edge cost $C(T)=\sum\nolimits_{e\in T} c(e)$ for a set $T\subseteq \mathcal{E}$. An arborescence is a directed, rooted tree in which all edges point away from the root; a minimal spanning forest for a digraph is the union of the arborescences which span the digraph, such that the total edge cost is as small as possible.  A matching for a graph $\mathcal{G}=(\mathcal{V},\mathcal{E})$ is a subset of edges in $\mathcal{E}$ which do not share common vertices. A maximum matching of $\mathcal{G}$ is a matching with the maximum number of edges among all possible matchings, whose size is called the matching number, denoted by $v(\mathcal{G})$. A minimum
cost maximum matching is the maximum matching with the edge cost as small as possible. 
Given a maximum matching $M^*$ of $\mathcal{B}(\mathcal{S}_1,\mathcal{S}_2,\mathcal{E}_{\mathcal{S}_1,\mathcal{S}_2})$, a vertex is said to be matched w.r.t. $M^*$, if it belongs to an edge in $M^*$, otherwise it is unmatched; a vertex is said to be right-matched (resp. left-matched) w.r.t. $M^*$, if it belongs to $\mathcal{S}_2$ (resp. $\mathcal{S}_1$) and in $V_R(M^*)$ (resp. $V_L(M^*)$); otherwise it is right-unmatched (resp. left-unmatched). We say $\mathcal{B}(\mathcal{S}_1,\mathcal{S}_2,\mathcal{E}_{\mathcal{S}_1,\mathcal{S}_2})$ has a perfect matching, 
 if there are no unmatched vertices w.r.t any maximum matching.

\subsection{Structural Controllability}
Consider a network system whose dynamic is captured by
\begin{equation}
\label{plant Eq}
x(t+1)=Ax(t)+Bu(t),
\end{equation}
where $x(t) \in \mathbb{R}^n$ is the state vector, $u(t) \in \mathbb{R}^q$ is the input vector, $A \in \mathbb{R}^{n\times n}$ and $B \in \mathbb{R}^{n \times q}$ are respectively the state transition matrix and input matrix. In practical, the exact values of entries of $A$ and $B$ might be hard to know.
Hence, let $\bar A \in {\{ 0,1\} ^{n \times n}}$  and $\bar B \in {\{ 0,1\} ^{n \times q}}$ be binary matrices representing the sparsity patterns of matrices $A$ and $B$, where $1$ denotes a free parameter and $0$ a zero entry. We call a matrix as {\emph{structured matrix}}, if every of its entry is either a fixed zero or a free parameter. For two structured matrices $\bar A_1$ and $\bar A_2$ with the same dimensions, we say $\bar A_1\subseteq \bar A_2$, if whenever $\bar A_{1{ij}}\ne 0$ implies $\bar A_{2ij}\ne 0$.


Let $\mathcal{X}$, $\mathcal{U}$ denote the sets of state vertices and input vertices respectively, i.e., $\mathcal{X}=\{x_1,...,x_n\}$, $\mathcal{U}=\{u_1,...,u_q\}$. Denote the edges by $\mathcal{E}_{\mathcal{X},\mathcal{X}}(\bar A)=\{(x_i,x_j): \bar A_{ji}\ne 0\}$,
$\mathcal{E}_{\mathcal{U},\mathcal{X}}(\bar B)=\{(u_j,x_i): \bar B_{ij}\ne 0\}$. An edge (i.e., link) $e$ is said to be \emph{state edge (link)} if $e \in \mathcal{X} \times \mathcal{X}$, and \emph{input edge (link)} if $e \in \mathcal{U}\times\mathcal{X}$. Let $\mathcal{D}(\bar A,\bar B)= (\mathcal{X}\cup \mathcal{U}, \mathcal{E}_{\mathcal{X},\mathcal{X}}(\bar A)\cup \mathcal{E}_{\mathcal{U},\mathcal{X}}(\bar B))$ be the system digraph associated with $(\bar A, \bar B)$; moreover, $\mathcal{D}(\bar A)=(\mathcal{X},\mathcal{E}_{\mathcal{X},\mathcal{X}})$, $\mathcal{D}(\bar B)=(\mathcal{X}\cup \mathcal{U}, \mathcal{E}_{\mathcal{U},\mathcal{X}})$. In the following, we sometimes simplify $\mathcal{E}_{\mathcal{X},\mathcal{X}}(\bar A)$ by $\mathcal{E}_{\mathcal{X},\mathcal{X}}$, $\mathcal{E}_{\mathcal{U},\mathcal{X}}(\bar B))$ by $\mathcal{E}_{\mathcal{U},\mathcal{X}}$, if no confusion is made. 

We say $(\bar A, \bar B)$ is structurally controllable if there exists a realization $(A,B)$ with the sparsity pattern of $(\bar A, \bar B)$ such that $(A,B)$ is controllable in the numerical sense.
 For the system $(\bar A, \bar B)$ in (\ref{plant Eq}) and its system digraph $\mathcal{D}(\bar A,\bar B)$, a state vertex $x\in\mathcal{X}$ is said to be input-reachable, if there exists at least one path from one of the input vertices $u\in\mathcal{U}$ to $x$ in $\mathcal{D}(\bar A,\bar B)$. Decomposing $\mathcal{D}(\bar A)=(\mathcal{X},\mathcal{E}_{\mathcal{X},\mathcal{X}})$ into SCCs, an SCC having no incoming edges from other SCCs to its vertices is called a source SCC. A source SCC is said to be a non input-reachable source SCC if none of its vertices is input-reachable. A stem is a path from an input vertex $u\in {\cal U}$ to a state vertex $x\in {\cal X}$. By connecting a stem and a collection of disjoint cycles with state or input links in the system digraph $\mathcal{D}(\bar A,\bar B)$, we get a cactus.

The generic rank of a structured matrix $M$ is the maximum rank $M$ can achieve as the function of its free parameters. We denote by ${\rm grank}(M)$ as the generic rank of $M$.

The following lemma characterizes the structural controllability, which can be found in \cite{Sergio_Pequito_2016}, \cite{multivariable_graph_1988} etc.


\begin{lemma}\label{Lemma 1}
Given a pair $(\bar A, \bar B)$, let $\mathcal{D}(\bar A,\bar B)= (\mathcal{X}\cup \mathcal{U}, \mathcal{E}_{\mathcal{X},\mathcal{U}}\cup \mathcal{E}_{\mathcal{U},\mathcal{X}})$ and $\mathcal{B}(\bar A,\bar B)= \mathcal{B}(\mathcal{X}\cup \mathcal{U}, \mathcal{X}, \mathcal{E}_{\mathcal{X},\mathcal{X}}\cup \mathcal{E}_{\mathcal{U},\mathcal{X}})$ be its system digraph and bipartite graph respectively. The following statements are equivalent: \begin{itemize}
\item[i).] The pair $(\bar A, \bar B)$ is structurally controllable;
\item[ii).]   $\mathcal{D}(\bar A,\bar B)$ can be spanned by a collection of disjoint cacti;
\item[iii).]  (a) every state vertex is input-reachable in $\mathcal{D}(\bar A,\bar B)$; \\~
              (b) there is a maximum matching for $\mathcal{B}(\bar A,\bar B)$, such that every state vertex is right-matched.
\item[iv).]  (a) there is a path from $\cal U$ to every $x\in {\cal X}$ in ${\cal D}(\bar A, \bar B)$;     \\~
             (b) ${\rm grank}([\bar A, \bar B])=n$.
             \end{itemize}
\end{lemma}

 {\emph{Structural perturbation.}} We call the addition or deletion of links or vertices to/from the digraph associated with a structured system, including state links (vertices) and input links (vertices), as structural perturbations. The vertex deletion removes a vertex as well as all links (edges) incident to or from such vertex from the original plant.

\subsection{Problem Statements}

Consider the structured system $(\bar A, \bar B)$ given in (\ref{plant Eq}). Without loss of generality, assume $||\bar B||_0>0$.  We are interested in the following problems related to structural perturbations subject to structural controllability:


{\bf Problem 1.} (Minimal cost link insertion) If $(\bar A, \bar B)$ is not structurally controllable,  determine the minimal cost of link set from a given link set (including state links and input links) whose insertion makes it structurally controllable.

{\bf Problem 2.} (Minimal cost link deletion) If $(\bar A, \bar B)$ is structurally controllable,  determine the minimal cost of link set (including state links and input links) whose deletion makes it structurally uncontrollable.

{\bf Problem 3.} (Minimal cost actuator deletion) If $(\bar A, \bar B)$ is structurally controllable, determine the minimal cost of actuator set whose deletion makes it structurally uncontrollable.

The mathematical formulations of the above problems can be found in the corresponding sections subsequently.  The heterogeneous costs imposed on different links or actuators are natural settings, noting that for practical networks, different links or actuators may incur different importance, difficulty or budgets to be added/deleted to/from a system. The above problems are inherently combinatorial optimization problems. Rigorous analysis to reveal their computation complexities is what we mainly pursue in this paper.

The above link insertion/deletion or actuator deletion problems can also be understood in the following way: regard each state variable as a follower, each input as a leader, and the non-zero entries in $\bar A$ and $\bar B$  as communication links among followers and from leaders to followers respectively, i.e., forming a leader-follower multi-agent system \cite{Multi-agent_2004}.  Then, the corresponding link intersetion/deletion or actuator deletion problems can be seen as adding/removing communication links or removing leaders in the associated multi-agent system.

\section{Minimum Cost Link Insertion Problem}

 In the minimal cost link insertion problem, given a pair $(\bar A, \bar B)$, and $(A_s,B_s)$ with the same dimensions as $(\bar A, \bar B)$ such that $\mathcal{E}_{{\cal{X}},{\cal{X}}}(\bar A_s), \mathcal{E}_{{\cal{U}},{\cal{X}}}(\bar B_s)$ denote the sets of candidate state edges and input edges that can be added to the original system,  respectively. Each candidate edge $e \in \mathcal{E}_{{\cal{X}},{\cal{X}}}(\bar A_s)\cup \mathcal{E}_{{\cal{U}},{\cal{X}}}(\bar B_s)$ is assigned a non-negative cost $c(e)\ge 0$.   We intend to select a subset of links with minimum cost from $\mathcal{E}_{{\cal{X}},{\cal{X}}}(\bar A_s)\cup \mathcal{E}_{{\cal{U}},{\cal{X}}}(\bar B_s)$, such that the resulting system is structurally controllable. This problem is formulated as
\[\begin{array}{l}
\mathop {\min }\limits_{\Delta \bar A \subseteq \bar A_s,\Delta \bar B  \subseteq  \bar B_s} {\kern 1pt} {\kern 1pt} {\kern 1pt} {\kern 1pt} {\kern 1pt} {\kern 1pt} {\kern 1pt} {\kern 1pt} {\kern 1pt} {\kern 1pt} {\kern 1pt} {\kern 1pt} {\kern 1pt} {\kern 1pt} {\kern 1pt} {\kern 1pt} {\kern 1pt} {\kern 1pt} {\kern 1pt} {\kern 1pt} {\kern 1pt} {\kern 1pt} {\kern 1pt} {\kern 1pt} {\kern 1pt} {\kern 1pt} {\kern 1pt} {\kern 1pt} {\kern 1pt} {\kern 1pt} {\kern 1pt}  { \sum\limits_{e \in {{\mathcal{E}}_{{\mathcal{X}},{\mathcal{X}}}}(\Delta \bar A) \cup {{\mathcal{E}}_{{\mathcal{U}},{\mathcal{X}}}}( \Delta\bar B)} {c(e)} } \\ {\tag{Problem 1}}
s.t.{\kern 1pt} {\kern 1pt} {\kern 1pt} {\kern 1pt} {\kern 1pt} {\kern 1pt} {\kern 1pt} {\kern 1pt} {\kern 1pt} {\kern 1pt} {\kern 1pt} {\kern 1pt} {\kern 1pt} {\kern 1pt} {\kern 1pt} {\kern 1pt} {\kern 1pt} {\kern 1pt} {\kern 1pt} {\kern 1pt} {\kern 1pt} {\kern 1pt} {\kern 1pt} (\bar A \vee \Delta \bar A,\bar B \vee \Delta \bar B) {\kern 2pt} \text{is structurally controllable}
\end{array}\]
where $\vee$ is the point-wise OR operation for binary matrices, i.e, ${(\bar M \vee \bar N)_{ij}}=\bar M_{ij} \vee \bar N_{ij}$.

For simplifying description, by adding the setting $c(e)=0$ for $e\in \mathcal{E}_{\mathcal{X},\mathcal{X}}(\bar A) \cup \mathcal{E}_{\mathcal{U},\mathcal{X}}(\bar B)$ (the rest weights remain the same), and denoting ${\mathcal{E}^{can}_{{\mathcal{X}},{\mathcal{X}}}}\triangleq \mathcal{E}_{\mathcal{X},\mathcal{X}}(\bar A_s) \cup \mathcal{E}_{\mathcal{X},\mathcal{X}}(\bar A) $ and  ${\mathcal{E}^{can}_{{\mathcal{U}},{\mathcal{X}}}}\triangleq  {\mathcal{E}_{{\mathcal{U}},{\mathcal{X}}}}(
\bar B)\cup {\mathcal{E}_{{\mathcal{U}},{\mathcal{X}}}}(\bar B_s)$, Problem 1 is equivalent to the following problem
\[\begin{array}{l}
\mathop {\min }\limits_{ {\mathcal{E}_{{\mathcal{X}},{\mathcal{X}}}}(\Delta \bar A) \subseteq {\mathcal{E}^{can}_{{\mathcal{X}},{\mathcal{X}}}},   {\mathcal{E}_{{\mathcal{U}},{\mathcal{X}}}}(\Delta \bar B)  \subseteq  {\mathcal{E}^{can}_{{\mathcal{U}},{\mathcal{X}}}}} {\kern 1pt} {\kern 1pt} {\kern 1pt} {\kern 1pt} {\kern 1pt} {\kern 1pt} {\kern 1pt} {\kern 1pt} {\kern 1pt} {\kern 1pt} {\kern 1pt} {\kern 1pt} {\kern 1pt} {\kern 1pt} {\kern 1pt} {\kern 1pt} {\kern 1pt} {\kern 1pt} {\kern 1pt} {\kern 1pt} {\kern 1pt} {\kern 1pt} {\kern 1pt} {\kern 1pt} {\kern 1pt} {\kern 1pt} {\kern 1pt} {\kern 1pt} {\kern 1pt} {\kern 1pt} {\kern 1pt}  { \sum\limits_{e \in {{\mathcal{E}}_{{\mathcal{X}},{\mathcal{X}}}}(\Delta \bar A) \cup {{\mathcal{E}}_{{\mathcal{U}},{\mathcal{X}}}}( \Delta\bar B)} {c(e)} } \\
s.t.{\kern 1pt} {\kern 1pt} {\kern 1pt} {\kern 1pt} {\kern 1pt} {\kern 1pt} {\kern 1pt} {\kern 1pt} {\kern 1pt} {\kern 1pt} {\kern 1pt} {\kern 1pt} {\kern 1pt} {\kern 1pt} {\kern 1pt} {\kern 1pt} {\kern 1pt} {\kern 1pt} {\kern 1pt} {\kern 1pt} {\kern 1pt} {\kern 1pt} {\kern 1pt} (\Delta \bar A, \Delta \bar B) {\kern 2pt} \text{is structurally controllable}
\end{array}.\]
Denote the above problem by $\mathcal{P}_{ins}^0({\mathcal{E}^{can}_{{\mathcal{X}},{\mathcal{X}}}},{\mathcal{E}^{can}_{{\mathcal{U}},{\mathcal{X}}}},C)$.
Moreover, by setting $c(e)=0$ for all $e\in \mathcal{E}^{can}_{\mathcal{X},\mathcal{X}}$,  Problem 1 collapses to the minimal cost input selection problems discussed in \cite{A_Olshevsky_2015}, \cite{Sergio_Pequito_2016}, \cite{Sergio_Pequito_2015_cost} under various cost $c(e)$ for $e\in \mathcal{E}^{can}_{\mathcal{U},\mathcal{X}}$. Those problems, as subproblems of Problem 1 where only input links can be inserted,  can be solved in polynomial time as shown in \cite{A_Olshevsky_2015}, \cite{Sergio_Pequito_2016}, \cite{Sergio_Pequito_2015_cost}. However, the following theorem reveals that Problem 1 is NP-hard in general. Such distinction is the essential
difference between the link insertion problem discussed in this
paper and the input selection problems in the existing literature.

\begin{theorem}\label{Theorem 0}
The minimal cost link insertion problem
(Problem 1) is NP-hard with identical link weights.
\end{theorem}

{\bf{Proof:}} We show a polynomial time reduction from the {\emph{Hamiltonian path problem}} to Problem 1.

\begin{figure}
  \centering
  \subfigure[]{
    \includegraphics[width=1.4in]{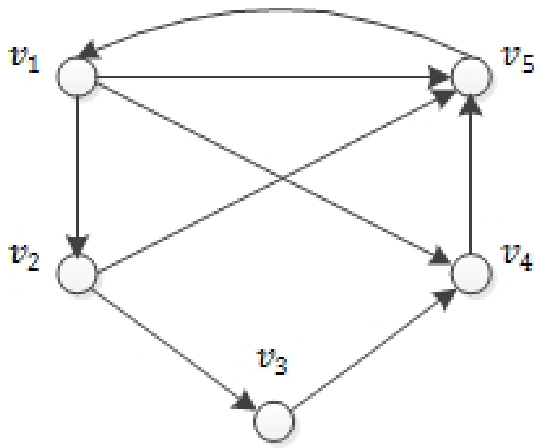}}
 \hspace{0.4in}
  \subfigure[]{
    \includegraphics[width=1.5in]{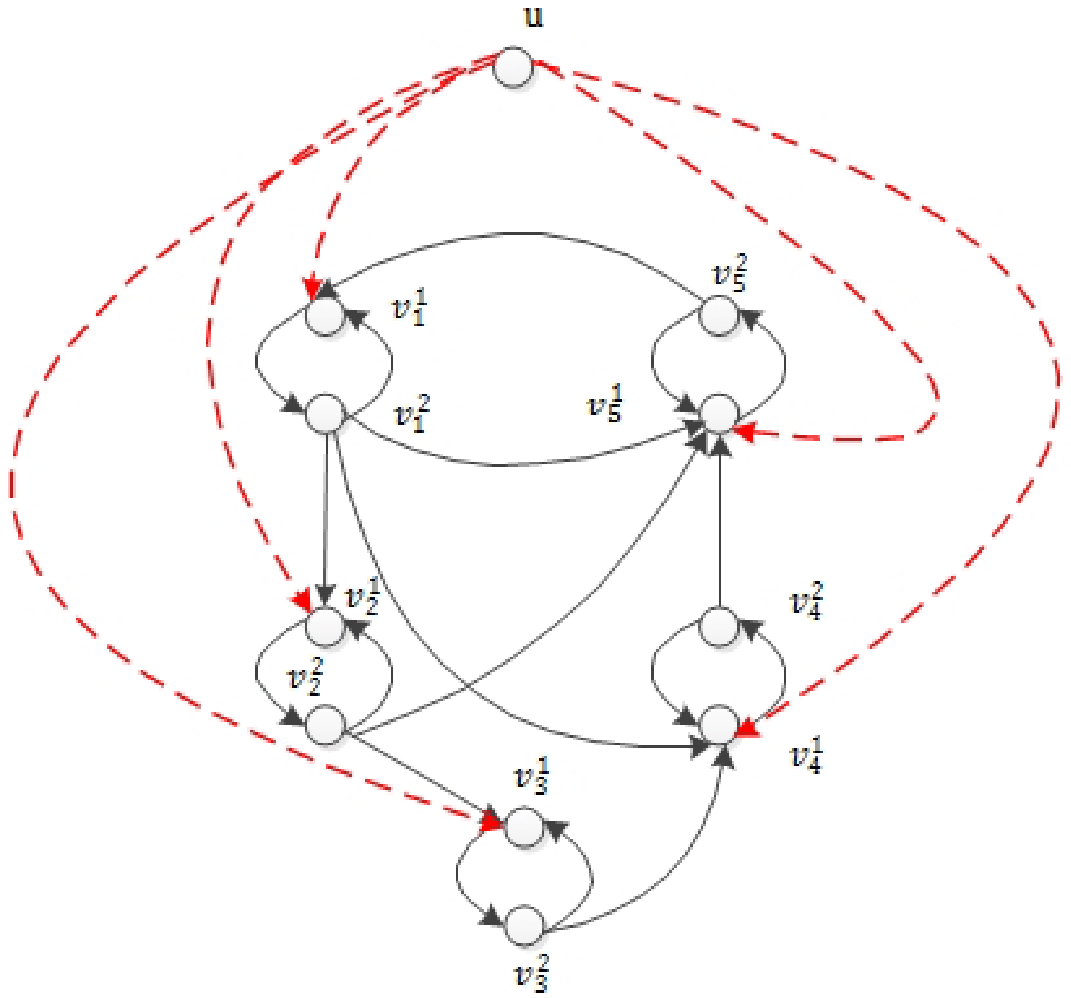}}
  \caption{\label{fignew1} Example of the construction from a digraph $\mathcal{G}$ to $\mathcal{P}_{ins}^0({\mathcal{E}}_p,{\mathcal{E}}_{{\mathcal{U}},{\mathcal{V}_p}},C)$. In Fig. \ref{fignew1},   (a) is the digraph $\mathcal{G}$, (b) is the auxiliary graph $\mathcal{G}_{p}=(\mathcal{V}_{p},{\mathcal{E}}_{p})$ (solid arrows) along with the input edge candidate ${\mathcal{E}}_{{\mathcal{U}},{\mathcal{V}_p}}$ (dotted arrows). }
  \label{fig:subfig} 
\end{figure}

A Hamiltonian path in a directed graph is a path visiting each vertex exactly once. Determining whether such paths exist in graphs is Hamiltonian path problem, which is NP-complete \cite{DB_West_graph}.
Now, given an arbitrary digraph $\mathcal{G}=(\mathcal{V},{\mathcal{E}})$, where $\mathcal{V}=\{v_1,...,v_n\}$, construct an auxiliary graph $\mathcal{G}_{p}=(\mathcal{V}_{p},{\mathcal{E}}_{p})$ where ${{\mathcal{V}_p}} = \{ v_1^1,v_1^2,...,v_n^1,v_n^2\}$, obtained from $\mathcal{G}$ by replacing each vertex $v_i$ of $\mathcal{G}$ with a cycle containing two vertices  $v_i^1$ and $v_i^2$, and letting $v_i^1$ have all the in neighbors as $v_i$ and $v_i^2$ all the out neighbors as $v_i$ for $i=1,...,n$. Add a single input vertex ${\mathcal{U}}=\{u\}$ to $\mathcal{G}_{p}$, connect $u$ to all $v_i^1$ vertices of $\mathcal{G}_{p}$ for $i=1,...,n$, and assign unit cost to each edge of the resulting digraph. Finally, map the obtained digraph to $\mathcal{P}_{ins}^0({\mathcal{E}}_p,{\mathcal{E}}_{{\mathcal{U}},{\mathcal{V}_p}},C)$, where ${{\mathcal{E}}_{{\mathcal{U}},{\mathcal{V}_p}}} = \{(u,v_i^1):i=1,...,n\}$, $C=\{c(e)=1: e \in {\mathcal{E}}_p \cup {{\mathcal{E}}_{{\mathcal{U}},{\mathcal{V}_p}}}\}$. See Fig. \ref{fignew1} for illustrating of such construction. {\footnote{The reason of duplicating each vertex of $\cal G$ is to make sure that the constructed $\mathcal{P}_{ins}^0({\mathcal{E}}_p,{\mathcal{E}}_{{\mathcal{U}},{\mathcal{V}_p}},C)$ is always feasible. If the reduced problem is not feasible (which can be verified in polynomial time), then the corresponding Problem 1 can be trivially solved in polynomial time.}}

It is easy to see that a structurally controllable system can be obtained from the edge candidates $({\mathcal{E}}_p,{\mathcal{E}}_{{\mathcal{U}},{\mathcal{V}_p}})$, as in any case the diagraph $({\mathcal{V}_p} \bigcup {\mathcal{U}}, {\mathcal{E}}_p\bigcup {\mathcal{E}}_{{\mathcal{U}},{\mathcal{V}_p}})$ itself can be covered by a cactus.
 We declaim that the optimal cost of $\mathcal{P}_{ins}^0({\mathcal{E}}_p,{\mathcal{E}}_{{\mathcal{U}},{\mathcal{V}_p}},C)$ is no more than $2n$, if and only if $\mathcal{G}$ has a Hamiltonian path.

In one direction, when $\mathcal{G}$ has a Hamiltonian path, denoted by $v_{p_1}\rightarrow ...\rightarrow v_{p_n}$, where $\{p_1,...,p_n\}$ is a perturbation of $\{1,...,n\}$,  there is a stem $u \to v_{{p_1}}^1 \to v_{{p_1}}^2 \to v_{{p_2}}^1... \to v_{{p_n}}^1 \to v_{{p_n}}^2$ with size $2n$ in the auxiliary graph $({\mathcal{V}_p} \bigcup {\mathcal{U}}, {\mathcal{E}}_p\bigcup {\mathcal{E}}_{{\mathcal{U}},{\mathcal{V}_p}})$, which corresponds to a structurally controllable system with total cost $2n$.

In the other direction, suppose there exists a structurally controllable system with total cost no more than $2n$ obtained from the edge collection $({\mathcal{E}}_p,{\mathcal{E}}_{{\mathcal{U}},{\mathcal{V}_p}})$, denoted by $\Sigma^*$. 
  Noting that there are $2n$ state vertices in $\Sigma^*$, at least $2n$ edges are needed to make every state vertex input reachable. In addition, applying condition ii) of
  Lemma \ref{Lemma 1}  to $\Sigma^*$, it follows that no cycle can exist in $\Sigma^*$. That is because, if there exists a cycle, then at least one state vertex has in-degree at least $2$ (one is from the cycle, the other is from the cactus), which leads to a total cost no less than $2+2n-1=2n+1$. Hence, the diagraph associated with $\Sigma^*$ must be spanned by a stem, denoted by $u \to v_{{s_1}}^1 \to v_{{s_1}}^2 \to v_{{s_2}}^1 \to... \to v_{{s_n}}^1 \to v_{{s_n}}^2$, where $\{s_1,...,s_n\}$ is a perturbation of $\{1,...,n\}$. Consequently, $v_{s_1} \to v_{s_2} \to...\to v_{s_n}$ forms a Hamiltonian path of $\mathcal{G}$. 

 Since the above reduction can be implemented in polynomial time, and determining whether there exists a Hamiltonian path in graph $\mathcal{G}$ is NP-complete, it concludes that verifying whether the optimal cost of $\mathcal{P}_{ins}^0({\mathcal{E}}_p,{\mathcal{E}}_{{\mathcal{U}},{\mathcal{V}_p}},C)$ is no more than $2n$ is NP-complete too. Therefore, Problem 1 is NP-hard. $\hfill\blacksquare$

 Following a similar argument of the proof of Theorem \ref{Theorem 0}, it leads to the following corollary.

 \begin{corollary}\label{Corollary 1} In Problem 1, provided the input topology $\mathcal{E}_{\mathcal{U},\mathcal{X}}(\Delta \bar B)$ is prescribed, it is NP-hard to determine $\mathcal{E}_{\mathcal{X},\mathcal{X}}(\Delta \bar A)\subseteq \mathcal{E}^{can}_{\mathcal{X},\mathcal{X}}$ with the minimal cost such that $(\Delta \bar A, \Delta \bar B)$ is structurally controllable; or equivalently, if $c(e)=0$ for all $e\in \mathcal{E}^{can}_{\mathcal{U},\mathcal{X}}$, Problem 1 is NP-hard.
  \end{corollary}

 {\emph{Proof:}} The proof is a slight modification of the proof of Theorem \ref{Theorem 0}. Given an arbitrary digraph $\cal{G}$ with vertex number $n$, construct the same auxiliary graph $({\mathcal{V}_p} \bigcup {\mathcal{U}}, {\mathcal{E}}_p\bigcup {\mathcal{E}}_{{\mathcal{U}},{\mathcal{V}_p}})$ with ${{\cal{U}}=\{u\}}$ in the same way as the proof of Theorem \ref{Theorem 0}, and add an extra vertex $z$ to the auxiliary graph along with an edge $(z,u)$. Denote the obtained graph by ${\cal{G}}_p'$. The difference is that now we regard $z$ as the only input vertex, while the rest vertices ${\cal{V}}_p \bigcup \{u\}$ as state vertices. Let the input link $(z,u)$ be fixed, i.e., setting $w((z,u))=0$, and all the rest edges of the auxiliary graph have unit cost. In such a construction, the corresponding minimal cost link insertion problem is always feasible, as the new auxiliary graph ${\cal{G}}_p'$ itself is always structurally controllable.
 Next, following similar analysis to the proof of Theorem \ref{Theorem 0}, it holds that there exists a structurally controllable system with total cost no more than $2n$ for ${\cal{P}}^0_{ins}$ associated with ${\cal{G}}_p'$, if and only if $\cal{G}$ has a Hamiltonian path. The latter problem is NP-complete. Hence, the result of Corollary  \ref{Corollary 1} follows immediately.  $\hfill\blacksquare$

Theorem \ref{Theorem 0} and Corollary \ref{Corollary 1} make it clear that determining the minimum cost structurally controllable network topology from a given collection of links is NP-hard, and it is still NP-hard to do so when the input topology is prescribed, even with identical link costs. These intractability results are in sharp contrast to the minimal input selection problems (in terms of the total cost of input links, see \cite{A_Olshevsky_2015}, \cite{Sergio_Pequito_2016}, \cite{Sergio_Pequito_2015_cost}) for a fixed autonomous network topology, which can be solved in polynomial time. Such distinction may result form the fact that, the latter problems are computationally equivalent to the corresponding {\emph {maximum matching}} problems for bipartite graphs as suggested in \cite{complexity_minimum_input}, while the former problem is not easier than the Hamiltonian path problem as revealed in the proof of Theorem 1. 
After a deeper insight, it seems that such distinction might result from the admission of adding state links, noticing that state edge addition involves both the start vertex and the end one of an edge with regard to the connectivity or the matching properties (Lemma \ref{Lemma 1}), while adding input edge merely needs to consider the status of the end state vertex of the added edge.

As for approximation, there is a $2$-approximation algorithm for Problem 1, which is a natural combination of the minimal spanning forest (arborescence) algorithm and the minimum cost maximum matching algorithm; see Algorithm \ref{alg0} and Theorem \ref{Theorem add 1}. The basic idea of Algorithm \ref{alg0} is to find the minimal cost of additional edges to form a maximum matching to match all state vertices, on the basis of a minimal spanning forest, then eliminate some redundant edges which don't destruct the input-reachability of all state vertices.

\begin{algorithm} 
  {{{\small{
\caption{: Approximation algorithm for Problem 1} 
\label{alg0} 
\begin{algorithmic}[1] 
\REQUIRE $({{\mathcal{E}}^{can}_{{\mathcal{X}},{\mathcal{X}}}},{{\mathcal{E}}^{can}_{{\mathcal{U}},{\mathcal{X}}}},C)$ 
\ENSURE   Approximation solution to $\mathcal{P}_{ins}^0({{\mathcal{E}}^{can}_{{\mathcal{X}},{\mathcal{X}}}},{{\mathcal{E}}^{can}_{{\mathcal{U}},{\mathcal{X}}}},C)$  
\STATE    Determine the minimal spanning forest of digraph $({\mathcal{U}}\bigcup {{\mathcal{X}},{{\mathcal{E}}^{can}_{{\mathcal{X}},{\mathcal{X}}}}\bigcup {{{\mathcal{E}}^{can}_{{\mathcal{U}},{\mathcal{X}}}}}})$ rooted in $\mathcal{U}$ with edge cost $C$ such that there is no isolated state vertex, denoted by $\mathcal{{{T}}}$;
 \STATE    Let $C'\leftarrow C$, and update $C'$ by letting $c'(e)=0$, $\forall e \in E({\mathcal{T}})$; construct the bipartite graph $({\mathcal{U}}\bigcup {\mathcal{X}}, {\mathcal{X}}, {{\mathcal{E}}^{can}_{{\mathcal{X}},{\mathcal{X}}}}\bigcup {{{\mathcal{E}}^{can}_{{\mathcal{U}},{\mathcal{X}}}}})$ with edge cost in $C'$, and determine its minimum cost maximum matching such that every state vertex is right-matched, denoted by $\mathcal{M}$;
\STATE Let $C''\leftarrow C$, and update $C''$ by letting $c''(e)=0$, $\forall e \in \mathcal{M}$;  find the minimal spanning forest of digraph $({\mathcal{U}}\bigcup {\mathcal{X}}, E(\mathcal{T})\cup \mathcal{M})$ rooted in $\mathcal{U}$ with edge cost in $C''$, given by $\mathcal{T}'$;
\STATE Return the structured system $(\bar A,\bar B)$ with $\mathcal{D}(\bar A, \bar B)=(\mathcal{U}\cup \mathcal{X}, E(\mathcal{T}^{'})\cup \mathcal{M})$. 
\end{algorithmic}}}
}}
\end{algorithm}

 \begin{theorem}\label{Theorem add 1}
If Problem 1 is feasible, Algorithm \ref{alg0} is a $2$-approximation to Problem 1 with complexity $\mathcal{O}((|\mathcal{U}|+|\mathcal{X}|)^3)$.
\end{theorem}

{\bf{Proof:}} Let $\mathcal{G}_{opt}$ be the digraph associated with the optimal solution to Problem 1.  As every state vertex is input-reachability, there must exist a spanning forest in $\mathcal{G}_{opt}$ with no isolated state vertices, given by $\mathcal{T}_{opt}$. By definition, $C(E(\mathcal{T}_{opt}))\ge C(E(\mathcal{T}))$. In addition, every state vertex should be right-matched by some maximum matching of the bipartite graph associated with $\mathcal{G}_{opt}$, denoted by $\mathcal{M}_{opt}$. Since
$\mathcal{M}_{opt} \subseteq E({{\mathcal{E}}^{can}_{{\mathcal{X}},{\mathcal{X}}}}\cup {{\mathcal{E}}^{can}_{{\mathcal{U}},{\mathcal{X}}}})$, and every edge with cost set $C'$ has a cost not larger than that of the corresponding edge with cost set $C$, it is clear that $C(\mathcal{M}_{opt}) \ge C'(\mathcal{M})=\sum\nolimits_{e \in \mathcal{M}}{c'(e)}$. Noticing that $E(\mathcal{T}')\subseteq E(\mathcal{T})\cup \mathcal{M}$, it follows $C(E(\mathcal{T}')\cup \mathcal{M})\le C(E(\mathcal{T})\cup \mathcal{M})=C(E(\mathcal{T}))+ C'(\mathcal{M}) \le C(E(\mathcal{T}_{opt}))+ C(\mathcal{M}_{opt}) \le 2 C(E(\mathcal{G}_{opt}))$. Hence, Algorithm \ref{alg0} achieves a $2$-approximation to Problem 1.


As for computation complexity, Steps 1 and 3 can be implemented using Edmonds' algorithm  in time $\mathcal{O}((|\mathcal{U}|+|\mathcal{X}|)|{{\mathcal{E}}^{can}_{{\mathcal{X}},{\mathcal{X}}}}\cup {{\mathcal{E}}^{can}_{{\mathcal{U}},{\mathcal{X}}}}|)$ \cite{DB_West_graph}. Step 2 costs $\mathcal{O}((|\mathcal{U}|+|\mathcal{X}|)^3)$ complexity using Hungarian algorithm \cite{DB_West_graph}. The rest steps have linear complexity. To sum up, Algorithm \ref{alg0} incurs in $\mathcal{O}((|\mathcal{U}|+|\mathcal{X}|)^3)$.  $\hfill\blacksquare$

\begin{figure}
\begin{minipage}[t]{1\linewidth}
\centering
\includegraphics[width=3.3in]{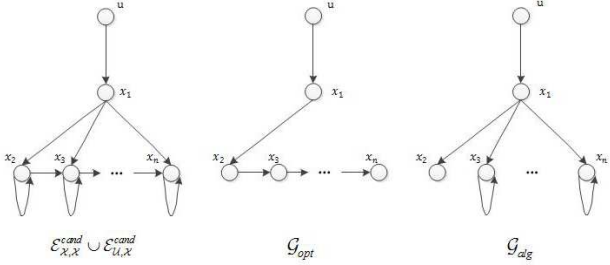}
\caption{\label{fignew2} Example for the worst case performance of Algorithm \ref{alg0}.}
\label{fig:side:a}
\end{minipage}%
\end{figure}

The bound `$2$' in Theorem \ref{Theorem add 1} is tight, which we can see by the example illustrated by Fig. \ref{fignew2},
in which every edge has unit cost. The optimal solution $\mathcal{G}_{opt}$ has a total cost $n$, while Algorithm \ref{alg0} might select a solution like $\mathcal{G}_{alg}$, whose cost is $2n-2$. As a consequence, the approximation factor ${\rho _n} = \frac{{2n - 2}}{n} = 2 - \frac{2}{n}$, which leads to $\mathop {\lim }\limits_{n \to \infty } {\rho _n} = 2$.
It is not difficult to see that, under each of the following two scenes, Algorithm \ref{alg0} always returns the optimal solution: (i). every state vertex $x\in {\cal X}$ has a zero-cost self-loop; (ii). $({\cal X},\mathcal{E}^{can}_{\mathcal{X},\mathcal{X}})$ can be covered by a strongly-connected subgraph with zero cost.

\begin{remark}
 (Iterative improvement of Algorithm \ref{alg0}) A natural direction to improve Algorithm \ref{alg0} is to implement Algorithm $1$ iteratively, and in each iteration, perturb an edge of the previous obtained $\mathcal{T}'$ to reconstruct a new spanning forest $\mathcal{T}$ (the cardinality of potential edges is at most $(|\mathcal{U}|+|\mathcal{X}|)|\mathcal{X}|$), and pick the edge with the largest decrease in the return value of Algorithm \ref{alg0} (at the price of increasing computation burden). It is easy to see that the optimal solution of the example of Fig. \ref{fignew2} can be obtained through such iterative improvement. However, such implement does not guarantee to return an optimal solution but may encounter suboptimal solutions.
 \end{remark}

Even though the optimal link insertion problem with controllability constraint is in general NP-hard, as shown in Theorem \ref{Theorem 0}, there are some restricted cases under which Problem 1 has polynomial time complexity. Particulary, under the scenario where there is no restriction on the insertable links and each link has $0-1$ cost, i.e., $\mathcal{E}^{can}_{\mathcal{X},\mathcal{X}}={\cal X}\times {\cal X}$, $\mathcal{E}^{can}_{\mathcal{U},\mathcal{X}}={\cal U}\times {\cal X}$, and $c(e)=0$ or $1$, $e\in \mathcal{E}^{can}_{\mathcal{X},\mathcal{X}}\cup \mathcal{E}^{can}_{\mathcal{U},\mathcal{X}}$,  Problem 1 can be solved in polynomial time. An equivalent formulation of the aforementioned problem is given as follows
\[\begin{array}{l}
{\rm Given } (\bar A, \bar B), ||\bar B||_0>0, {\rm determine}\\
\mathop {\min }\limits_{\Delta \bar A \in {{\{ 0, 1\} }^{n \times n}},\Delta \bar B  \in {{\{ 0, 1\} }^{n \times q}}} {\kern 1pt} {\kern 1pt} {\kern 1pt} {\kern 1pt} {\kern 1pt} {\kern 1pt} {\kern 1pt} {\kern 1pt} {\kern 1pt} {\kern 1pt} {\kern 1pt} {\kern 1pt} {\kern 1pt} {\kern 1pt} {\kern 1pt} {\kern 1pt} {\kern 1pt} {\kern 1pt} {\kern 1pt} {\kern 1pt} {\kern 1pt} {\kern 1pt} {\kern 1pt} {\kern 1pt} {\kern 1pt} {\kern 1pt} {\kern 1pt} {\kern 1pt} {\kern 1pt} {\kern 1pt} {\kern 1pt} {\left\| {\Delta \bar A} \right\|_0}{\rm{ + }}{\left\| {\Delta \bar B } \right\|_0}\\
s.t.{\kern 1pt} {\kern 1pt} {\kern 1pt} {\kern 1pt} {\kern 1pt} {\kern 1pt} {\kern 1pt} {\kern 1pt} {\kern 1pt} {\kern 1pt} {\kern 1pt} {\kern 1pt} {\kern 1pt} {\kern 1pt} {\kern 1pt} {\kern 1pt} {\kern 1pt} {\kern 1pt} {\kern 1pt} {\kern 1pt} {\kern 1pt} {\kern 1pt} {\kern 1pt} (\bar A \vee \Delta \bar A,\bar B \vee \Delta \bar B) {\kern 2pt} \text{is structurally controllable}
\end{array}\]
where ${\left\| { M } \right\|_0}$ denotes the zero norm.
For details, see the conference paper \cite{Y_Zhang_2017}. A similar result is also independently obtained in \cite{same} recently.

\section{Minimum cost link deletion problem}
In the above section we have considered the link addition to system (\ref{plant Eq}). Now we consider the link deletion from system (\ref{plant Eq}).  Given $(\bar A, \bar B)$ in (\ref{plant Eq}),  each link $e\in \mathcal{E}_{\mathcal{X},\mathcal{X}}(\bar A)\cup \mathcal{E}_{\mathcal{U},\mathcal{X}}(\bar B)$ has a non-negative link cost $c(e)\ge 0$. The {\emph{minimal cost link deletion problem }} aims to minimize the cost of the set of links whose removal from $\mathcal{E}_{\mathcal{X},\mathcal{X}}(\bar A)\cup \mathcal{E}_{\mathcal{U},\mathcal{X}}(\bar B)$ precludes the existence of a structurally controllable system constructed from the rest links. Formally, it can be formulated as
\[\begin{array}{l}
\mathop {\min }\limits_{\Delta \bar A \subseteq \bar A,\Delta \bar B  \subseteq  \bar B} {\kern 1pt} {\kern 1pt} {\kern 1pt} {\kern 1pt} {\kern 1pt} {\kern 1pt} {\kern 1pt} {\kern 1pt} {\kern 1pt} {\kern 1pt} {\kern 1pt} {\kern 1pt} {\kern 1pt} {\kern 1pt} {\kern 1pt} {\kern 1pt} {\kern 1pt} {\kern 1pt} {\kern 1pt} {\kern 1pt} {\kern 1pt} {\kern 1pt} {\kern 1pt} {\kern 1pt} {\kern 1pt} {\kern 1pt} {\kern 1pt} {\kern 1pt} {\kern 1pt} {\kern 1pt} {\kern 1pt}  { \sum\limits_{e \in {{\mathcal{E}}_{{\mathcal{X}},{\mathcal{X}}}}(\Delta \bar A) \cup {{\mathcal{E}}_{{\mathcal{U}},{\mathcal{X}}}}( \Delta\bar B)} {c(e)} } \\ {\tag{Problem 2}}
s.t.{\kern 1pt} {\kern 1pt} {\kern 1pt} {\kern 1pt} {\kern 1pt} {\kern 1pt} {\kern 1pt} {\kern 1pt} {\kern 1pt} {\kern 1pt} {\kern 1pt} {\kern 1pt} {\kern 1pt} {\kern 1pt} {\kern 1pt} {\kern 1pt} {\kern 1pt} {\kern 1pt} {\kern 1pt} {\kern 1pt} {\kern 1pt} {\kern 1pt} {\kern 1pt} (\bar A \backslash \Delta \bar A,\bar B \backslash \Delta \bar B) {\kern 2pt} \text{is structurally uncontrollable}
\end{array},\]
where for two binary matrices $\bar M$ and $\bar N$, $\backslash$ is the entry-wise subtraction operation, satisfying ${(\bar M{\rm{ \backslash }}\bar N)_{ij}} = {\rm{1}}$ if and only if ${\bar M_{ij}} = {\rm{1}}$, ${\bar N_{ij}} = 0$. 

Intuitively, the solution to Problem 2 measures how hard it is to destroy the structural controllability of system $(\bar A, \bar B)$. Any link failures with a total cost less than the optimum of Problem 2 can't destruct the system controllability. In this sense, the optimum of Problem  is a measure of robustness against link failures/deletions w.r.t controllability.


\begin{theorem}\label{Theorem 3}
The minimal cost link deletion problem (Problem 2) is NP-hard even with identical link costs. 
\end{theorem}

Before presenting the proof, let us we analyze the general solution to Problem 2 with each link having unit weight.  According to Lemma \ref{Lemma 1}, it is straightforward to see that, given a pair $(\bar A, \bar B)$,  the minimum number of edges whose deletion destroys structural controllability is equal to the minimum number of edges whose deletion destroys the input-reachability of $\mathcal{D}(\bar A, \bar B)=(\mathcal{X}\cup \mathcal{U}, \mathcal{E}_{\mathcal{X},\mathcal{X}}\cup \mathcal{E}_{\mathcal{U},\mathcal{X}})$ or the maximum matching of $\mathcal{B}(\bar A,\bar B)= \mathcal{B}(\mathcal{X}\cup \mathcal{U}, \mathcal{X}, \mathcal{E}_{\mathcal{X},\mathcal{X}}\cup \mathcal{E}_{\mathcal{U},\mathcal{X}})$. For further discussion, the following notions related to the graph connectivity and matching are needed. Readers can refer to \cite{M_C_Dour_2015}, \cite{DB_West_graph}, \cite{R_Zenkusen_2009} for more details.

In the following, let $\mathcal{D}=(\mathcal{V},\mathcal{E})$ be a digraph with $s$ and $t \in \mathcal{V}$ being the source and the sink of $\mathcal{D}$, and every edge $(u,v) \in \mathcal{E}$ mapping to a capacity $c_{uv}>0$.

\begin{definition} (Minimum cut) Given the digraph $\mathcal{D}$, an $s\text{-}t$ cut is a set of edges whose removal leads to the non-existence of paths from $s$ to $t$. The minimum cut problem is to determine an $s\text{-}t$ cut with the minimal sum of edge capacities. \end{definition}

%


Introduce a virtual source $\bar u$ to $\mathcal{D}(\bar A, \bar B)$, such that there is an edge from $\bar u$ to every $u\in \mathcal{U}$, i.e., $ \mathcal{E}_{\bar u,\mathcal{U}}=\{(\bar u, u_i): i=1,...,q\}$, denoted the resulting digraph by $\mathcal{D}(\bar A, \bar B, \bar u)=(\mathcal{X}\cup \mathcal{U}\cup \{\bar u\}, \mathcal{E}_{\mathcal{X},\mathcal{X}}\cup \mathcal{E}_{\mathcal{U},\mathcal{X}} \cup \mathcal{E}_{\bar u, \mathcal{U}})$.  Assign the capacity as :
$c_{e}=1$ if $e\in \mathcal{E}_{\mathcal{X},\mathcal{X}}\cup \mathcal{E}_{\mathcal{U},\mathcal{X}}$, and $c_{\bar u u_i}=\sum\nolimits_{\{ x:(u_i,x) \in \mathcal{E}_{\mathcal{U},\mathcal{X}}\}}{{c_{u_ix}}}$.
Then, it is clear that, for a given $x_i \in \mathcal{X}$, the minimum edges whose deletion destroys the input-reachability of $x_i$ equals to the minimum $\bar u\text{-}x_i$ cut in $\mathcal{D}(\bar A, \bar B, \bar u)$, denoted by $\lambda (\bar u,{x_i})$.
Let ${T_{cut}}(\bar A,\bar B)$ be the minimum number of edges whose deletion  destroys the input-reachability of at least one state vertex in $\mathcal{D}(\bar A, \bar B)$. According to the above, it is easy to see that
$${T_{cut}}(\bar A,\bar B) = \mathop {\min }\nolimits_{1 \le i \le n} {\kern 1pt} \lambda (\bar u,{x_i}).$$

\begin{definition} (\cite{R_Zenkusen_2009}) ($1$-blocker, minimum cost $1$-blocker, and matching preclusion) Given an undirected graph $\mathcal{G}=(\mathcal{V},\mathcal{E})$ with matching number $v(\mathcal{G})$, a subset $\mathcal{E}_s \subseteq \mathcal{E}$ is a $d$-blockers of $\mathcal{G}$ if it satisfies $v(\mathcal{V}, \mathcal{E}\backslash \mathcal{E}_s)\leq v(\mathcal{G})-d$. If each edge in $\mathcal{E}$ has a non-negative cost, the $1$-blocker with the minimum cost is the minimum cost $1$-blocker among all possible 1-blockers. Specifically, when $d=1$ and $\mathcal{G}$ has a perfect matching, the minimum edge size of $1$-blocker of $\mathcal{G}$ is also called the matching preclusion number (see \cite{M_C_Dour_2015} for another definition).  \end{definition}  


The following lemma characterizes the NP-completeness of the $1$-blocker problem and the matching preclusion number.

\begin{lemma}[ Theorem 3.3 of \cite{R_Zenkusen_2009};\cite{M_C_Dour_2015}]\label{lemma 6}
For a bipartite graph $\mathcal{B}(\mathcal{S}_1,\mathcal{S}_2,\mathcal{E}_{\mathcal{S}_1,\mathcal{S}_2})$ and a given integer $r$, it is NP-complete to determine whether there exists a $1$-blocker of size at most $k$; when $|\mathcal{S}_1|=|\mathcal{S}_2|$, it is NP-complete to decide whether the matching preclusion number of $\mathcal{B}(\mathcal{S}_1,\mathcal{S}_2,\mathcal{E}_{\mathcal{S}_1,\mathcal{S}_2})$ is at most $r$.
\end{lemma}

For a structurally controllable pair $(\bar A, \bar B)$, it is clear that the minimum number of edges whose deletion destroys the matching condition is equal to the minimum $1$-blocker of the bipartite graph $\mathcal{B}(\bar A,\bar B)$, denoted by $T_{bl}({\mathcal B}(\bar A,\bar B))$. Let $d_{\bar c}(\bar A, \bar B)$ be the minimum edges whose deletion destroys the structurally controllability of $(\bar A, \bar B)$. Then, it follows that
\begin{equation}
\label{min_deletion}
{d_{\bar c}}(\bar A,\bar B) = \min \left\{ {{T_{cut}}(\bar A,\bar B),T_{bl}({\mathcal B}(\bar A,\bar B))} \right\}.
\end{equation}

From the above, $T_{cut}(\bar A,\bar B)$ can be determined in polynomial time by solving $|\mathcal{X}|$ max-flow problems in $\mathcal{D}(\bar A, \bar B, \bar u)$ according to the well-known Max-flow min-cut theorem \cite{DB_West_graph}, more specifically, with complexity of $|\mathcal{X}|^2(|\mathcal{E}_{\mathcal{X},\mathcal{X}}|+|\mathcal{E}_{\mathcal{U},\mathcal{X}}|)$ using the Edmonds-Karp algorithm  \cite{DB_West_graph}. By Lemma \ref{lemma 6}, the minimum $1$-blocker problem is NP-hard in general. However, we can not conclude that Problem 2 is NP-hard yet. That's because, the resulting bipartite graph $\mathcal{B}(\bar A,\bar B)$ has some inherent structure, such that we can not declaim that determining $T_{bl}({\mathcal B}(\bar A,\bar B))$ is NP-hard. In particular, $\mathcal{B}(\bar A,\bar B)$ corresponds to a digraph $\mathcal{D}(\bar A, \bar B)$ where every vertex $x_i\in \mathcal{X}$ is reachable from at least one $u_j \in \mathcal{U}$. What is more, even if it is NP-hard to determine $T_{bl}({\mathcal B}(\bar A,\bar B))$, we have to verify whether its value is less than ${T_{cut}}(\bar A,\bar B)$, whose size usually varies with $\mathcal{D}(\bar A,\bar B)$ but not being constant. {\emph{The difficulty is therefore to construct a transformation
from the $1$-blocker problem of general bipartite graphs to an instance of Problem 2, while exploring an explicit relationship of size between the minimum cut and the minimum $1$-blocker involved therein.}} In the following we provide a rigorous proof satisfying the above requirements.


{\bf{Proof of Theorem \ref{Theorem 3}:}}  Given a structurally controllable pair $(\bar A, \bar B)$ and an integer $r$, for arbitrary pairs $(\bar A_s, \bar B_s)$ with feasible dimensions, it can be verified whether ${\left\| {{{\bar A}_{\rm{s}}}} \right\|_0} + {\left\| {{{\bar B}_{\rm{s}}}} \right\|_0} \leq r$  and $(\bar A \backslash \bar A_s, \bar B \backslash \bar B_s)$ is structurally controllable in polynomial time. Therefore, the decision version of Problem 2 is NP.

To prove the NP-hardness, we build an instance of Problem 2 starting from the matching preclusion number problem of a generic bipartite graph. Let $\mathcal{B}(\mathcal{S}_1,\mathcal{S}_2,\mathcal{E}_{\mathcal{S}_1,\mathcal{S}_2})$ be bipartite with a perfect matching and $|\mathcal{S}_1|=|\mathcal{S}_2|=n$. Construct a structured system $(\bar A, \bar B)$ as: the state vertex set $\mathcal{X}=\{x_1,...,x_n\}$, the input vertex set $\mathcal{U}=\{u_1,...,u_n\}$, and $\mathcal{E}_{\mathcal{U},\mathcal{X}}=\{(u_i,x_j):(s_i,s_j)\in \mathcal{E}_{\mathcal{S}_1,\mathcal{S}_2}\}$, $\mathcal{E}_{\mathcal{X},\mathcal{X}}=\emptyset$. That is, the corresponding $\bar A$, $\bar B$ are respectively
$$\bar A = {0_{n \times n}},{\bar B_{ij}} = \left\{ \begin{array}{l}
1,{\kern 1pt} {\kern 1pt} {\kern 1pt} if{\kern 1pt} {\kern 1pt} {\kern 1pt} ({s_j},{s_i}) \in {\mathcal{E}_{\mathcal{S}_1,\mathcal{S}_2}}\\
0,{\kern 1pt} {\kern 1pt} {\kern 1pt} else
\end{array} \right.$$
Let $\hat {\mathcal B}(\bar A,\bar B)=\mathcal{B}(\mathcal{U}, \mathcal{X},\mathcal{E}_{\mathcal{U},\mathcal{X}})$. It is easy to see that the resulting system $(\bar A, \bar B)$ satisfies: ~\\
(i) every $x\in \mathcal{X}$ can be matched as $\mathcal{B}(\mathcal{S}_1,\mathcal{S}_2,\mathcal{E}_{\mathcal{S}_1,\mathcal{S}_2})$ has a perfect matching (so is with $\hat {\mathcal B}(\bar A,\bar B)$; ~\\
(ii) every $x\in \mathcal{X}$ is input-reachable, as every $x\in \mathcal{X}$ is matched by a $u\in \mathcal{U}$ w.r.t any perfect matching of $\hat {\mathcal B}(\bar A,\bar B)$. ~\\
Consequently, $(\bar A, \bar B)$ is structurally controllable.

According to the max-flow min-cut theorem and the structure property of digraph $(\mathcal{X}\cup\mathcal{U}, \mathcal{E}_{\mathcal{U},\mathcal{X}})$,  ${T_{cut}}(\bar A,\bar B) = \mathop {\min }\nolimits_{1 \le i \le n} {\kern 1pt} {\kern 1pt} {\kern 1pt} \deg ({x_i})$, where $\deg(x_i)$ denotes the in-degree of $x_i\in \mathcal{X}$, i.e., $\deg ({x_i}) = \sum\nolimits_{j = 1}^n {{{\bar B}_{ij}}} $.
From the property of matching preclusion number, it is valid that
\[T_{bl}(\hat {\mathcal B}(\bar A,\bar B)) \le \mathop {\min }\nolimits_{1 \le i \le n} {\kern 1pt} {\kern 1pt} {\kern 1pt} \deg ({x_i}){\rm{ = }}{T_{cut}}(\bar A,\bar B).\]
The left-hand relation is obvious as deleting all the input edges of an arbitrary vertex will certainly destroy a perfect matching. Then, according to (\ref{min_deletion}), $${d_{\bar c}}(\bar A,\bar B)\!\! =\!\! \min \left\{ {{T_{cut}}(\bar A,\bar B),T_{bl}(\hat {\mathcal B}(\bar A,\bar B))} \right\}\! =\! T_{bl}(\hat {\mathcal B}(\bar A,\bar B)).$$
Consequently, the minimum edge deletion to transform $(\bar A, \bar B)$ to be structurally uncontrollable is less than a given integer $r$, if and only if the matching preclusion number of the bipartite graph $\mathcal{B}(\mathcal{S}_1,\mathcal{S}_2,\mathcal{E}_{\mathcal{S}_1,\mathcal{S}_2})$ is below $r$. Since the latter is NP-complete, and the reduction can be done in polynomial time, it concludes that the decision version of Problem 2 is NP-complete, or alternatively, Problem 2 is NP-hard. This completes the proof. $\hfill\blacksquare$

From the above analysis, for approximation of Problem 2, we have the following conclusions.
\begin{theorem}\label{theorem 6}
If there exists a multiplicative factor $f(n)$ approximation algorithm for the minimal cost $1$-blocker problem, there is a $f(n)$-approximation algorithm for Problem 2, where $n$ is the input size of the corresponding problem.
\end{theorem}

{\bf{Proof:}}  Let each edge in ${{\mathcal{E}}_{{\mathcal{X}},{\mathcal{X}}}}(\bar A)\cup {{\mathcal{E}}_{{\mathcal{U}},{\mathcal{X}}}}(\bar B)$ have multiple costs (i.e, capacities). Following a similar argument to the analysis of unit link costs, denote the associated cost of minimum cut by ${{T^c_{cut}}(\bar A,\bar B)}$, which can be obtained in polynomially time, and the corresponding minimum cost $1$-block by $T^c_{bl}({\mathcal B}(\bar A,\bar B))$. Then, it can be seen that the optimum to Problem 2 is ${d^{opt}_{\bar c}}(\bar A,\bar B) = \min \left\{ {{T^c_{cut}}(\bar A,\bar B),T^c_{bl}({\mathcal B}(\bar A,\bar B))} \right\}$. If there is a $f(n)$-approximation algorithm for the minimal cost $1$-blocker problem and implementing such algorithm on $T^c_{bl}({\mathcal B}(\bar A,\bar B))$ returns $\bar T^c_{bl}({\mathcal B}(\bar A,\bar B))$, construct an algorithm which returns ${d^{alg}_{\bar c}}(\bar A,\bar B) = \min \left\{ {{T^c_{cut}}(\bar A,\bar B), \bar T^c_{bl}({\mathcal B}(\bar A,\bar B))} \right\}$. By definition, $\bar T^c_{bl}({\mathcal B}(\bar A,\bar B))\le f(n){T^c_{bl}}({\mathcal B}(\bar A,\bar B))$. Hence, it follows
\[\begin{array}{l}
{d^{alg}_{\bar c}}(\bar A,\bar B) \le \min \{{T^c_{cut}}(\bar A,\bar B),f(n){T^c_{bl}}({\mathcal B}(\bar A,\bar B))\}\\
 \le \min \{f(n){T^c_{cut}}(\bar A,\bar B),f(n){T^c_{bl}}({\mathcal B}(\bar A,\bar B))\} = f(n){d^{opt}_{\bar c}}(\bar A,\bar B)
\end{array}\] This finishes the proof. $\hfill\blacksquare$

In the reduction of the proof of Theorem \ref{Theorem 3}, since $\left\|\bar A\right\|_0=0$, the edges that can be deleted happen to be restricted in the input edges (i.e., $\mathcal{D}(\bar B)$), which immediately leads to the following corollary.

\begin{corollary}\label{corollary_2}
It is NP-hard to determine the minimum number of input links whose deletion destructs the structural controllability of a system.
\end{corollary}

\begin{remark}\label{remark 2} Corollary \ref{corollary_2} answers the hardness of determining the largest number of communication link (i.e., input link) failures a multi-agent system can robustly admit before structural controllability is preserved in \cite{M_A_Rahimian_failures_2013}. 
Theorem \ref{theorem 6} makes it clear that Problem 2 generally has the same multiplicative approximation factor as that of the minimal cost $1$-blocker problem. Readers can refer to \cite{R_Zenkusen_2010} for discussions on the latter problem.
\end{remark}  

\section{Minimal cost actuator deletion problem}
The former two sections have focused on the link addition/deletion to/from system (\ref{plant Eq}). In this section we consider the actuator deletion from system (\ref{plant Eq}). For $(\bar A, \bar B)$ in (\ref{plant Eq}), $q\ge 1$, let $\bar B_J$ be the submatrix of $\bar B$ formed by column vectors indexed by $J\subseteq \{1,...,q\}$.  Each input has a cost $c(i): \mathbb{N} \to \mathbb{R}_{\ge 0}$, measuring the importance of such input to the network, or the difficulty to be removed. Let $S=\{1,...,q\}$.

The minimal cost actuator deletion problem intends to determine the minimal cost of actuators whose deletion destroys structural controllability of the network.
This problem can be formulated as:
\[\begin{array}{l}
\mathop {\min {\kern 1pt} {\kern 1pt} {\kern 1pt} }\limits_{{J_c} \subseteq S} {\kern 1pt} \sum\limits_{i \in {J_c}} {c(i)} \\ \tag{Problem 3}
{\rm{s}}{\rm{.t}}{\rm{.}}{\kern 1pt} {\kern 1pt} {\kern 1pt} {\kern 1pt} (\bar A,{\bar B_{S\backslash {J_c}}}){\kern 1pt} {\rm{\quad is \quad structurally \quad uncontrollable}}
\end{array}\]

\begin{theorem} \label{NP_hard_input_deletion}
Problem 3 is NP-hard in the strong sense, even when each input actuates only one state vertex.
\end{theorem}

Notice that strong NP-hardness (NP-hard in the strong sense) implies that (unless P=NP) there cannot exist a fully polynomial-time
approximation scheme (FPTAS), i.e., an algorithm that solves a minimization problem within a factor of $1+\varepsilon$ of the optimal value in polynomial time of the input size and $1/\varepsilon$.  A problem is said to be strongly NP-complete, if it remains
so even when all of its numerical parameters are bounded by a polynomial in the length of the
input. A problem is said to be strongly NP-hard if a strongly NP-complete problem has a
polynomial reduction to it \cite{1985_spark}.  To prove the NP-hardness, some notion is introduced. The grith of a structured matrix $M$ is the minimal number of linearly dependent columns of $M$  \cite{1985_spark} (for a numerical matrix, the corresponding concept is called {\emph {spark}}).
 In the following proof, an {\emph{input removal set}} for $(\bar A, \bar B)$ is the set of inputs whose deletion causes structurally uncontrollability of the system. Denote $M_{J_1,J_2}$ as the submatrix of matrix $M$ formed by rows indexed by $J_1$ and columns indexed by $J_2$.


{\bf{Proof of Theorem \ref{NP_hard_input_deletion}:}} We adopt a reduction from the strongly NP-complete $k$-clique problem to an instance of Problem 3. A $k$-clique in a graph is a subgraph with any two of its vertexes being adjacent. The $k$-clique problem is to determine whether a undirected graph has a clique with size $k$.
Let ${\cal G}=({\cal V}, {\cal E})$ be a undirected graph, and $|{\cal V}|=n$, $|{\cal E}|=m$.   Denote the incidence matrix of ${\cal G}$ by $In({\cal G})$. Without loss of generality, assume that $\cal G$ is connected, and let $k> 4$ and satisfy \footnote{Combining the subsequent derivation,  Inequality (\ref{equality}) ensures that $k$ has the possibility to be the size of a clique of $\cal G$, noting that $\cal G$ is connected. Inequality (\ref{equality}) can be validated in polynomial time.} \begin{equation} \label{equality} {\small{ \left( {\begin{array}{*{20}{c}}
k\\
2
\end{array}} \right)}}+n-k\le m\end{equation} Construct an $(m+1)\times (m+1)$ structured matrix ${\cal C}({\cal G})$ as
{\small{
\begin{equation}{\cal C} (\cal G) \!\!=\!\! \left[ {\begin{array}{*{20}{c}}
{In(G)}&\vline& {{0_{n \times 1}}}\\
\hline
{{0_{(m{\rm{ + }}2{\rm{ + }}k - n - {\tiny{ \left( {\begin{array}{*{20}{c}}
k\\
2
\end{array}} \right)}}) \times m}}}&\vline& {{0_{(m{\rm{ + }}2{\rm{ + }}k - n - {\tiny{ \left( {\begin{array}{*{20}{c}}
k\\
2
\end{array}} \right)}}) \times 1}}}\\
\hline
{{1_{({\tiny{ \left( {\begin{array}{*{20}{c}}
k\\
2
\end{array}} \right)}} - k - 1) \times 1}}}&\vline& {\begin{array}{*{20}{c}}
{{0_{( {\tiny{\left({\begin{array}{*{20}{c}}
k\\
2
\end{array}} \right)}} - k - 2) \times 1}}}\\
{{1_{1 \times 1}}}
\end{array}}
\end{array}} \right]\end{equation}
}}As $k> 4$, it can be validated that ${\tiny{\left({\begin{array}{*{20}{c}}
k\\
2
\end{array}} \right)}} - k - 1> 0$, $m{\rm{ + }}2{\rm{ + }}k - n - {\tiny{ \left( {\begin{array}{*{20}{c}}
k\\
2
\end{array}} \right)}}>0$ from (\ref{equality}). Thus the construction is physically reasonable and ${\cal C}(\cal G)$ is square.

Construct an instance of Problem 3 as $\bar A = {{\cal C}}({\cal G})^\intercal$, $\bar B = {I_{m + 1}}$, with input costs
$$c(i) = \left\{ \begin{array}{l}
1,i = 1,...,m\\
m + 1,i = m + 1
\end{array} \right.$$

Obviously $(\bar A, \bar B)$ is structurally controllable. We declare that {\emph{the minimal cost input removal set for $(\bar A, \bar B)$ equals ${\tiny{{\left( \begin{array}{l}
k\\
2
\end{array} \right)}}}$, if and only if ${\cal G}$ has a $k$-clique.}}

To show this, an important property of the submatrix ${{\cal C}_1}({\cal G}) \triangleq {\cal C}{({\cal G})_{{C_1} \cup {C_3},S}}$ is utilized, with $C_1=\{1,...,n\}$, $C_3=\{m+k-{\tiny{ \left( {\begin{array}{*{20}{c}}
k\\
2
\end{array}} \right)}} +3,...,m+1\}$:  it demonstrates in \citep[Page 53]{1985_spark} that, matrix ${{\cal C}_1}({\cal G})$ has a girth with size ${\tiny{{\left( \begin{array}{l}
k\\
2
\end{array} \right)}}}$, if and only if ${\cal G}$ has a $k$-clique. 

For the one direction,  suppose that the minimal cost input removal set of $(\bar A, I)$ equals {${\tiny{{\left(\begin{array}{l}
k\\
2
\end{array} \right)}}}$. Denote the corresponding column index set by $J_c$. As ${\tiny{{\left( \begin{array}{l}
k\\
2
\end{array} \right)}}} \le m$, we have that $J_c\subseteq S$, and $m+1 \notin J_c$; otherwise $J_c$ has a cost no less than ${\tiny{{\left( \begin{array}{l}
k\\
2
\end{array} \right)}}}+1$. Let $S_{+}=\{1,...,m+1\}$.  Note that state vertices $\{{x_1},...,{x_m}\}$ are all out-neighbors of vertex $x_{m+1}$
from ${\cal D}(\bar A, \bar B)$. Hence, in the obtained system $(\bar A, \bar B_{S_{+}\backslash J_c})$ after removing inputs indexed by $J_c$, every state vertex is input-reachable. According to Lemma \ref{Lemma 1}, under such case, $(\bar A, \bar B_{S_{+}\backslash J_c})$ is structurally uncontrollable, if and only if
\begin{equation} \label{rank_deficient} {\rm grank}([\bar A, \bar B_{S_{+}\backslash J_c} ])<m+1.\end{equation}
Notice that every column of $\bar B_{S_{+}\backslash J_c}$ has only one nonzero entry. Therefore, (\ref{rank_deficient}) is equivalent to that
\begin{equation} \label{rank_second} {\rm grank}(\bar A_{J_c, S})<|J_c|,\end{equation}
Condition (\ref{rank_second}) also means that,  any $J_c \subseteq S$ making columns of ${({\bar A^\intercal})_{{C_1} \cup {C_3},{J_c}}}$ linearly dependent, is an input removal set for $(\bar A, \bar B)$.  As ${\tiny{{\left(\begin{array}{l}
k\\
2
\end{array} \right)}}}$ is the minimal cost input removal set, and each input in $S$ has unit cost, we have that ${\tiny{{\left(\begin{array}{l}
k\\
2
\end{array} \right)}}}$ is the girth of matrix ${({\bar A^\intercal})_{{C_1} \cup {C_3},S}}$ by definition, i.e., ${{\cal C}_1}({\cal G})$. According to the property of ${{\cal C}_1}({\cal G})$, this immediately leads to that graph $\cal G$ has a $k$-clique.

For the other direction, suppose  there is a $k$-clique in $\cal G$. By the property of ${{\cal C}_1}({\cal G})$, ${{\cal C}_1}({\cal G})$ has a girth ${\tiny\left( \begin{array}{l}
k\\
2
\end{array} \right)}$. Denote the column index set of such spark by $J_c\subseteq S$. It indicates that column vectors of ${({\bar A^\intercal})_{{C_1} \cup {C_3},{J_c}}}$ are linearly dependent, which immediately follows that ${\rm grank} ({\bar A_{{J_c},S}}) < |J_c|$. As a result, we have that
$${\rm grank}([\bar A, \bar B_{S_{+}\backslash J_c}])<m+1.$$
The above inequality leads to the uncontrollability of $(\bar A, \bar B_{S_{+}\backslash J_c})$. That is, $J_c$ is an input removal set with cost ${\tiny\left( \begin{array}{l}
k\\
2
\end{array} \right)}$. Moreover, as $|J_c|$ is the grith of ${{\cal C}_1}({\cal G})$, no other input removal set $J'_c\subseteq S$ with $|J'_c|<|J_c|$ exists.  On the other hand,
 notice that ${\tiny\left( \begin{array}{l}
k\\
2
\end{array} \right)}< m+1$. Hence, any input removal set containing $m+1$ will have a cost larger than $|J_c|$. Consequently, $J_c$ is the minimal input removal set with cost ${\tiny\left( \begin{array}{l}
k\\
2
\end{array} \right)}$.

The above reduction is within polynomial time. Combining the fact that $k$-clique problem is strongly NP-hard, the result follows.
$\hfill\blacksquare$

Several important remarks about the above proof should be noted here:

(i). While it has been proved that girth of a structured matrix  is NP-hard, the NP-hardness of Problem 3 can not be obtained directly from that fact. We explain it as follows.
For a given $(\bar A, I_n)$ with unit input cost, $\bar A\in \{0,1\}^{n\times n}$,  suppose ${\cal D}(\bar A)$ has $l$ source SCCs, denoted by ${\cal N}_1,...,{\cal N}_l$.
Following (\ref{rank_second}), the minimal cost input removal set equals the smaller value between $|J_0|$ (, $J_0\leftarrow \arg \min\nolimits_{J_c\subseteq S} {\rm grank}(\bar A^\intercal_{J_c})<|J_c|$) and $\min \nolimits_{i=1,...,l}|{\cal N}_i|$.  Even though the former value is NP-hard to determine, no explicit relation in size between the two aforementioned values can be found for a general matrix $\bar A$ as far as we know. In fact, there are many situations where the latter value is less than the former value, such that Problem 3 has polynomial time complexity (e.g., see Corollary {\ref{corollary_4}}). 

(ii). In the proof of Theorem {\ref{NP_hard_input_deletion}}, the matrix ${\cal C}_1({\cal G})$ originated from the construction of \cite{1985_spark}, where the author constructed  ${\cal C}_1({\cal G})$ to prove that girth is NP-hard. Here, we add a new column to ${\cal C}_1({\cal G})$ (i.e., the $(m+1)$-th column of ${\cal C}({\cal G})$), along with $m+2+k-n-{\tiny\left( \begin{array}{l}
k\\
2
\end{array} \right)}$ zero rows, and then assign specific costs to the inputs. With such specifically constructed ${\cal C}({\cal G})$, we demonstrate that the minimal cost input removal set equals the minimal size of cliques of $\cal G$ (rather than the value related to the SCCs of ${\mathcal D}(\bar A)$ mentioned in (i)).

(iii). If ${\cal C}_1(\cal G)$ is replaced with a general structured matrix, it can not be guaranteed that the corresponding statements in the proof of Theorem {\ref{NP_hard_input_deletion}} still hold.  Because of these reasons (also  (i) and (ii)),  as far as we know the explicit construction of ${\cal C}(\cal G)$ in the proof of Theorem {\ref{NP_hard_input_deletion}} is inevitable.

(iv). It is worthwhile to mention that Theorem \ref{NP_hard_input_deletion} does not indicate that determining the minimal number of actuators whose failure causes structural uncontrollability is NP-hard. Such problem is left for our further investigation.


%

The NP-hardness of Problem 3 does not rule out the possibility that under some restricted cases Problem 3 can be solved in polynomial time. We end this section by discussing one of such cases.
Suppose in a network system $(\bar A, \bar B)$, every state vertex has a self-loop, which is usually satisfied by physical systems \cite{NaturePhysics}, \cite{T_H_Summers_2016}. Suppose there are $l$ source SCCs in ${\cal D}(\bar A)$. For each source SCC ${\cal N}_i$, denote by $N({\cal N}_i)$ the neighbors of ${\cal N}_i$ in ${\cal D}(\bar A, \bar B)$ (hence $N({\cal N}_i)\subseteq {\cal U}$). Define $c^u_{i}=\sum \nolimits_{j\in \{j:u_j\in N({\cal N}_i)\}}c(j)$; that is, $c^u_{i}$ is the sum of costs of the inputs that are reachable to ${\cal N}_i$. We have the following conclusion.



\begin{corollary}\label{corollary_4}
Consider Problem 3 for a network with every state vertex having a self-loop. Problem 3 can be solved with complexity ${\cal O}(|{\cal X}|^2)$, and the optimal value equals $\min\nolimits_{1\le i\le l} c^u_{i}$.
\end{corollary}

{\bf Proof:} Since every state vertex has a self-loop, $[\bar A,\bar B_{S\backslash J_c}]$ is of full row generic rank $\forall J_c\subseteq S$. According to Lemma \ref{Lemma 1}, the minimal cost removal set equals the minimal cost of inputs whose deletion makes at least one state vertex input-unreachable. By the definition of source SCC, it suffices to see that such value is equal to $\min\nolimits_{1\le i\le l} c^u_{i}$. The complexity is dominated by the SCC-decomposition, which has complexity ${\cal O}(|\cal X|+|\cal E_{\cal X, \cal X}|)$, i.e., at most ${\cal O}(|{\cal X}|^2)$. $\hfill\blacksquare$

\section{Conclusions}
This paper addresses the problems of adding links with the minimal cost from a given set of links, including state links and input links, to make a network structurally controllable, and of removing links/actuators with the minimal cost to make a network structurally uncontrollable. We prove the NP-hardness of these problems. We also provide some approximation results for these related problems. The intractable results imply that it is generally hard to measure the `nearest distance' between structural controllability and structural uncontrollability in terms of number of links.   These results may serve an answer to the general hardness of optimally designing (modifying) a structurally
controllable network topology and of measuring controllability
robustness against link/actuator failures. Further work includes exploring more polynomial time algorithms to approximate these problems or determine optimal solutions to some of their subproblems.

\end{document}